\documentclass[11pt]{article}
\usepackage{epsfig}
\usepackage{subfigure}
\usepackage{authblk}
\usepackage{setspace}
\usepackage{amsmath}
\usepackage{float}
\usepackage{breqn}
\begin{document}
%\FLD{1}{6}{00}{28}{00}

 % Define commands to assure consistent treatment throughout document
 \newcommand{\eqnref}[1]{(\ref{#1})}
 \newcommand{\class}[1]{\texttt{#1}}
 \newcommand{\package}[1]{\texttt{#1}}
 \newcommand{\file}[1]{\texttt{#1}}
 \newcommand{\BibTeX}{\textsc{Bib}\TeX}
 \doublespacing

\title{\bf \vskip -1.5in Numerical Anisotropy in Finite Differencing}
\author{Adrian Sescu\thanks{ Department of Aerospace Engineering, Mississippi State University, 330 Walker at Hardy Rd, Mississippi State, MS 39762; {\it email}: sescu@ae.msstate.edu}\hspace{1.8mm}}
\affil{Department of Aerospace Engineering, Mississippi State University, MS 39762}

\date{}

\maketitle

\begin{abstract}
Numerical solutions to hyperbolic partial differential equations, involving wave propagations in one direction, are subject to several specific errors, such as numerical dispersion, dissipation or aliasing. In multi-dimensions, where the waves propagate in all directions, there is an additional specific error resulting from the discretization of spatial derivatives along grid lines. Specifically, waves or wave packets in multi-dimensions propagate at different phase or group velocities, respectively, along different directions. A commonly used term for the aforementioned multidimensional discretization error is the numerical anisotropy or isotropy error. In this review, the numerical anisotropy is briefly described in the context of the wave equation in multi-dimensions. Then, several important studies that were focused on optimizations of finite difference schemes with the objective of reducing the numerical anisotropy are discussed.

%\keywords{Finite differencing \and Numerical Anisotropy \and Partial differential equations \and Wave Propagation}
\end{abstract}

%\tableofcontents

\section{Introduction}
\label{intro}

Numerical anisotropy is a discretization error that is specific to numerical approximations of multidimensional hyperbolic partial differential equations (PDE). This error is often neglected, and the focus is directed toward the reduction of other types of discretization errors, such as numerical dissipation, dispersion or aliasing (e.g., Lele \cite{lele}, Tam and Webb \cite{tam1}, Kim and Lee \cite{Kim2}, Zingg et al. \cite{Zingg2}, Mahesh \cite{Mahesh}, Hixon \cite{hixon1}, Ashcroft and Zhang \cite{Ashcroft}, Fauconnier et al. \cite{Fauconnier} or Laizet and Lamballais \cite{Laizet}), or toward improving the accuracy of various time marching schemes (e.g., Hu et al. \cite{Hu}, Stanescu and Habashi \cite{Stanescu}, Mead and Renaut \cite{Mead}, Bogey and Bailly \cite{Bogey} or Berland et al. \cite{Berland}). There are several areas, however, where the numerical anisotropy can significantly affect the numerical solution based on finite difference or finite volume schemes (example include computational acoustics, computational electromagnetics, elasticity or seismology). The numerical anisotropy can be reduced by using, for example, one-dimensional high-resolution discretization schemes, multidimensional optimized difference schemes, or sufficiently fine grids. However, by increasing the number of grid points the computational time may increase considerably, while one-dimensional high-resolution difference schemes may generate spurious waves at the boundaries of the domain. Oftentimes, optimizations of multidimensional difference schemes are more effective.

High-order finite difference schemes that are optimized in one-dimension may not preserve their wavenumber resolution in multi-dimensional problems. These schemes may experience numerical anisotropy, because the dispersion characteristics along grid lines may not be the same as the dispersion characteristics associated with the diagonal directions. Over the years, several attempts to reduce the numerical anisotropy by various techniques were reported. A comprehensive analysis of the numerical anisotropy was performed in the book of Vichnevetsky~\cite{vich} where, among others, the two-dimensional wave equation was solved using two different finite difference schemes for the Laplacian operator. A considerable reduction of the numerical anisotropy was attained by weight averaging the two schemes. A slightly similar approach was previously used by Trefethen~\cite{trefethen} who used the leap frog scheme to solve the wave equation in two dimensions. Zingg and Lomax~\cite{Zingg1} performed optimizations of finite difference schemes applied to regular triangular grids, that give six neighbor points for a given node. They conducted comparisons between the newly derived schemes and conventional schemes that were discretized on square grids, and found that the numerical anisotropy can be significantly reduced by using triangular grids. Tam and Webb ~\cite{tam2} proposed an anisotropy correction to the finite difference representation of the Helmholtz equation. They derived an anisotropy correction factor using asymptotic solutions to the continuous equation and its finite difference approximation. 

Jo et al. \cite{Jo}, in the context of solving the acoustic wave equation, proposed a finite difference scheme over a stencil consisting of grid points from more than one direction, by linearly combining two discretizations of the second derivative operator. A notable reduction of the numerical anisotropy was obtined, but the numerical dispersion error was increased. Hustesdt et al. \cite{Hustesdt} proposed a two-staggered-grid finite difference schemes for the acoustic wave propagation in two dimensions, where the first derivative operator was discretized along the grid line and along the diagonal direction. Lin and Sheu~\cite{lin} explored the dispersion-relation-preserving concept of Tam and Webb~\cite{tam1} in two dimensions to optimize the first-order spatial derivative terms of a model equation that resembles the incompressible Navier-Stokes momentum equation. They approximated the derivative using a nine-point grid stencil, resulting in nine unknown coefficients. Eight of them were determined by employing Taylor series expansions, while the ninth one was determined by requiring that the two-dimensional numerical dispersion relation is the same as the exact dispersion relation. 

Kumar \cite{Kumar} derived isotropic finite difference schemes for the first and second derivatives in the context of symmetric dendritic solidification, and obtained a notable reduction of the numerical anisotropy. Patra and Kartunnen \cite{Patra} introduced several finite difference stencils for the Laplacian, Bilaplacian, and gradient of Laplacian, with the objective of improving the isotropic characteristics. Their stencils consisted of more grid points than the conventional schemes, but it was shown that the computational cost may decrease with more than 20$\%$ due to some gain in terms of stability. Stegeman et al. \cite{Stegeman} applied spectral analysis to evaluate the error in numerical group velocity (both the magnitude and the direction) of vorticity, entropy and acoustic waves, using the numerical solution to the linearized Euler equations in two-dimensions. They showed that a different measure of the group velocity error must be used to account for the error in the propagation direction of the waves. They also stressed that the numerical group velocity is more important than the numerical phase velocity in analyzing the errors associated with wave propagation. In a series of papers~\cite{sescu1,sescu2,sescu3,sescu4}, Sescu et al. proposed a technique to derive finite difference schemes in multi-dimensions with improved isotropy. The optimization performed in \cite{sescu1,sescu2,sescu3,sescu4} improved the isotropy of the wave propagation and, moreover, the stability restrictions of the multidimensional schemes in combination with either Runge-Kutta or linear multistep time marching methods were found to be more effective. They found that the stability restrictions are more favorable when using multidimensional schemes, even if they involve more grid points in the stencils. However, this was advantageous for low order schemes, such as those of second or fourth order of accuracy, but it was also shown that favorable stability restrictions can be obtained for higher order of accuracy schemes (sixth or eight) by increasing the isotropy corrector factor. The approach was extended to prefactored compact schemes by Sescu and Hixon \cite{sescu6,sescu7}. Beside reducing the numerical anisotropy, the new multidimensional compact schemes are computationally cheaper than the corresponding explicit multidimensional scheme defined on a same stencil.

In computational electromagnetics, there were many attempts to reduce the numerical anisotropy, by applying various techniques. Berini and Wu \cite{Berini} conducted a comprehensive analysis of the numerical dispersion and numerical anisotropy of finite difference schemes applied to transmission-line modeling (TLM) meshes. They found that, under certain circumstances, the time domain nodes introduce anisotropy into the dispersion characteristics of isotropic media, stressing the importance of developing schemes with improved isotropy. Gaitonde and Shang \cite{Gaitonde} proposed a class of high-order compact difference-based finite-volume schemes that minimizes the dispersion and isotropy error functions for the range of wavenumbers of interest. Sun and Trueman~\cite{Sun1} proposed an optimization of two-dimensional finite difference schemes, by considering additional nodes surrounding the point of differencing. They obtained a significant reduction in the numerical anisotropy, dispersion error and the accumulated phase errors over a broad bandwidth. Further optimizations of this scheme were performed in another paper of Sun and Trueman \cite{Sun2}. Koh et al. \cite{Koh} derived a two-dimensional finite-difference time-domain method, discretizing the Maxwell equations, to eliminate the numerical dispersion and anisotropy. They showed that the new algorithm has isotropic dispersion and resemble the exact phase velocity, whose isotropic property is superior to that of other existing schemes. Shen and Cangellaris~\cite{Shen} introduced a new stencil for the spatial discretization of Maxwell's equations. Compared to conventional second-order accurate FDTD scheme, their scheme experienced superior isotropy characteristics of the numerical phase velocity. They also showed that the Courant number cab be increased by using the newly derived schemes.  Kim et al. \cite{Kim} derived new three-dimensional isotropic dispersion-finite-difference time-domain schemes (ID-FDTD) based on a linear combination of the traditional central difference equation and a new difference equation using extra sampling points. Among all versions of the proposed finite-difference schemes, three of them showed improved isotropy of the wave propagation compared to the original scheme of the Yee \cite{Yee}. Kong and Chu \cite{Kong} introduced a new unconditionally-stable finite-difference time-domain method with low numerical anisotropy in three-dimensions. Compared with other finite-difference time-domain methods, the normalized numerical phase velocity of their proposed scheme was significantly improved, while the dispersion error and numerical anisotropy have been reduced.

This review will describe and discuss the numerical anisotropy in the framework of wave equation, and will present some of the most important optimizations of finite difference schemes in the context of reducing the numerical anisotropy. In section II, the dispersion error and the numerical anisotropy existing in finite difference discretizations of the wave equation are introduced and discussed. In section III, several approaches to reduce the numerical anisotropy, that were developed over the years by various research groups, are reviewed and discussed. Concluding remarks are included in section IV.

\section{Dispersion Error and Numerical Anisotropy}
\label{sec:1}

Let us consider the centered finite difference approximation of the spatial derivative, which contains both the explicit and the implicit (or compact) parts:

\begin{equation}\label{cvb}
\sum_{k=1}^{N_c}\alpha_k(u_{j+k}' + u_{j-k}') + u_{j}' =
\frac{1}{h}
\left(
\sum_{k=1}^{N_e}a_k(u_{j+k} - u_{j-k}) 
\right) +
O(h^n),
\end{equation}
where the gridfunctions are $u_j = u(x_j)$ for $1 \le j \le N$, the derivatives are denoted by a prime, $u_j'$, $h$ is the space step, and $\alpha_k$ and $a_k$ are given coefficients. If $N_{c} = 0$ the scheme is termed explicit, while compact schemes (also known as implicit or Pade schemes), by contrast, have $N_{c} \neq 0$ and require the solution of a matrix equation to determine the derivatives along a grid line. Conventionally, the coefficients $\alpha_k$ and $a_k$ are chosen to provide the largest possible exponent, $n$, in the truncation error, for a given stencil width, but in some instances some of these coefficients are determined to provide improved dispersion characteristics of the scheme. Table~\ref{t1} includes some of these weights for various explicit and compact finite difference schemes: explicit classical second order scheme (E2), explicit classical fourth order scheme (E4), explicit classical sixth order scheme (E6), dispersion-relation-preserving scheme of Tam and Webb~\cite{tam1}, compact classical fourth order scheme (C4), optimized tridiagonal compact scheme of Haras and Ta'asan~\cite{haras} (Haras), optimized pentadiagonal scheme of Lui and Lele~\cite{lui3} (Lui) and spectral-like pentadiagonal compact scheme of Lele~\cite{lele} (Lele). The prefactored compact scheme of Hixon~\cite{hixon1,hixon2} is also included here in the form:

\begin{eqnarray}\label{jjj}
a u_{j+1}^{F'} + c u_{j-1}^{F'} + (1-a-c) u_{j}^{F'} =
\frac{1}{h}
\left[
b u_{j+1} - (2b-1) u_j - (1-b) u_{j-1}) 
\right], \nonumber \\
c u_{j+1}^{B'} + a u_{j-1}^{B'} + (1-a-c) u_{j}^{B'} =
\frac{1}{h}
\left[
(1-b) u_{j+1} - (2b-1) u_j - b u_{j-1}) 
\right],
\end{eqnarray}
where $F$ and $B$ stand for 'forward' and 'backward', respectively (in a predictor-corrector time marching framework). For sixth order accuracy, $a=1/2-1/(2\sqrt{5})$, $b=1-1/(30 a)$ and $c=0$. The leading order term in the truncation error of a finite difference scheme depends on the choice of the coefficients and the $(n+1)$st derivative of the function $u$.

To study the wavenumber characteristics of finite difference schemes, consider a periodic domain in real space, $x \in [0,L]$, with $N$ uniformly spaced points(the spatial step size is $h=L/N$). The discrete Fourier transform of $u$ is given as
$
\hat{u}_m =
\frac{1}{N}\sum_{j=1}^{N} u_j e^{-ik_m x_j}$ with $m=-N/2,...,N/2-1,
$
where the wavenumber is $k_m = 2\pi m/L$. The $m$th component of the discrete Fourier transform of $u'$ denoted $\hat{u}_m'$ is simply $ik_m \hat{u}_m$. Taking the discrete Fourier transform of equation (\ref{cvb}) implies that

\begin{equation}\label{}
(\hat{u}_m')_{num} =
iK(k_m h)\hat{u}_m,
\end{equation}
where the numerical wavenumber is given as

\begin{equation}\label{}
K(z) =
\frac
{\sum_{n=1}^{N_e}2a_n \sin{(nz)}}
{1+\sum_{n=1}^{N_c}2\alpha_n \cos{(nz)}}.
\end{equation}

\begin{figure}[H]% order of placement preference: here, top, bottom
\centerline{
 \includegraphics[width=0.5\textwidth]{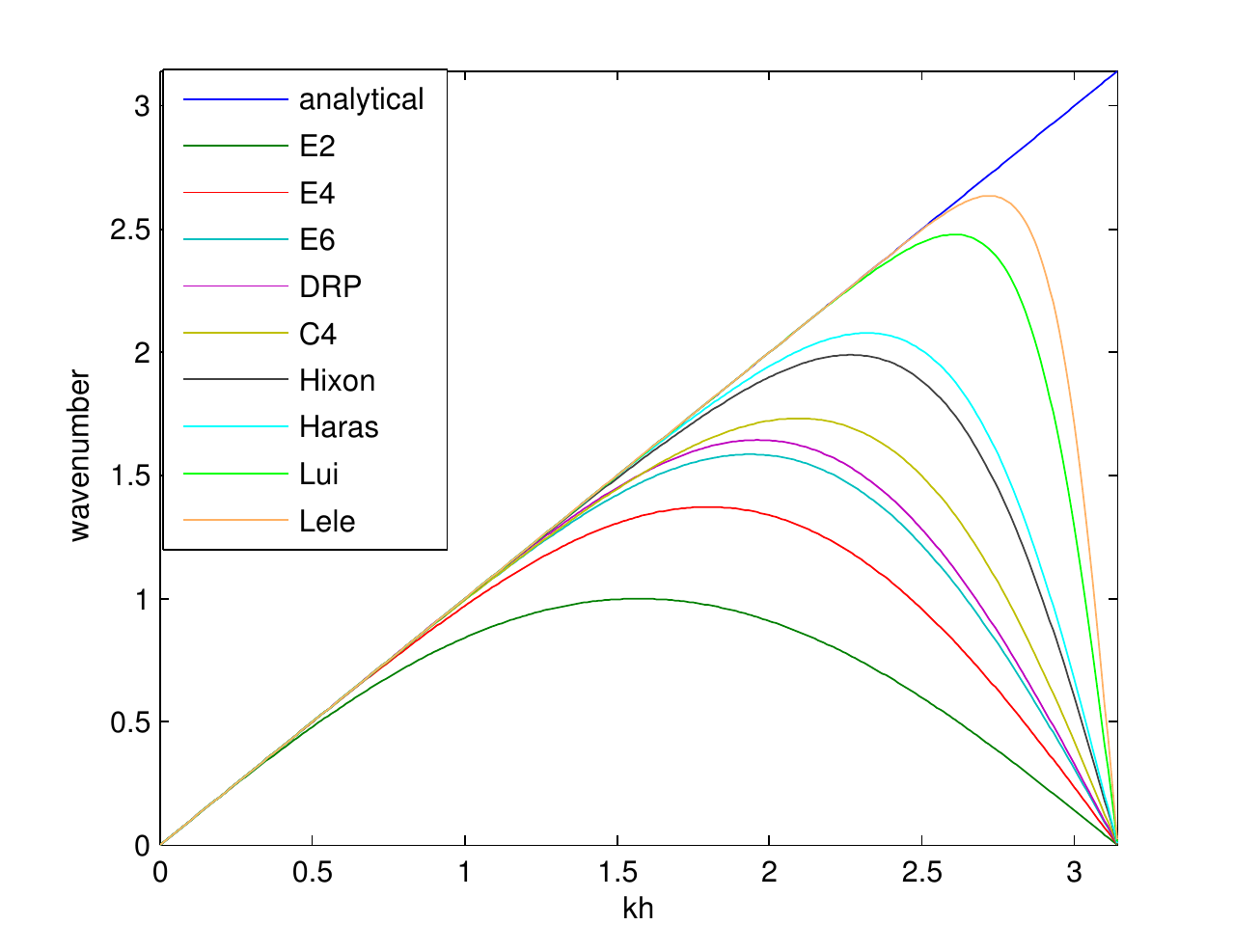}
           }
 \caption{Numerical wavenumber compared to the analytical wavenumber.}
 \label{f1}
\end{figure}

Figure \ref{f1} shows the numerical wavenumber for various explicit and compact schemes, corresponding to those given in table \ref{t1}. The numerical wavenumber is compared to the analytical wavenumber which is represented by the straight line in figure \ref{f1}. As one can notice, the compact schemes are superior to the explicit schemes; however, compact schemes are computationally more demanding because large matrices have to be inverted.

\begin{figure}[H]
  \begin{center}
    \mbox{
       \label{f2a}\includegraphics[width=0.35\textwidth]{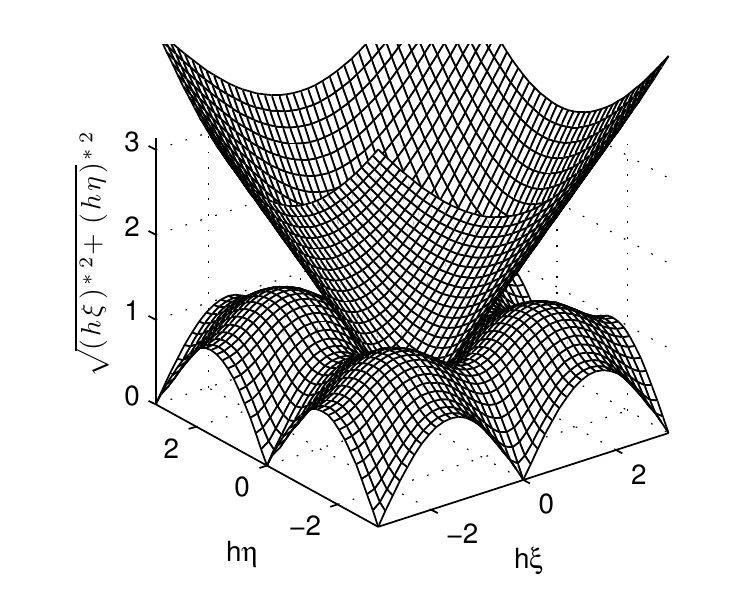}
       \label{f2b}\includegraphics[width=0.35\textwidth]{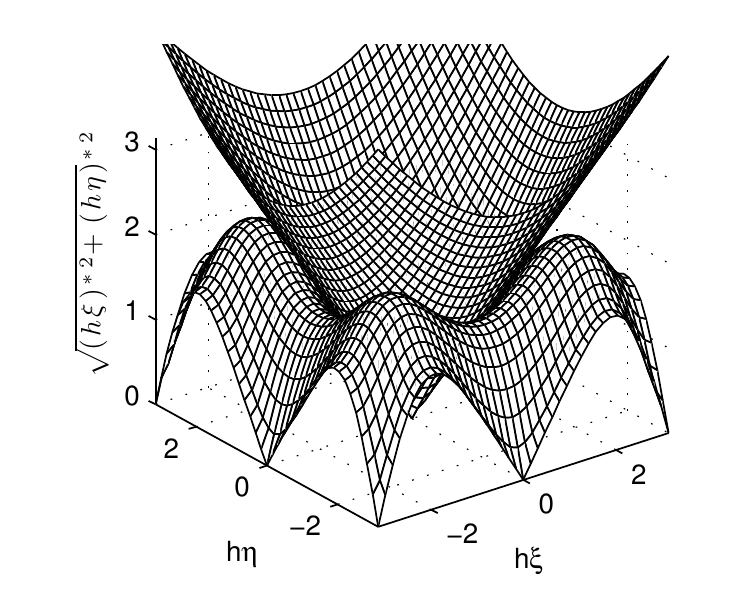}
       \label{f2b}\includegraphics[width=0.35\textwidth]{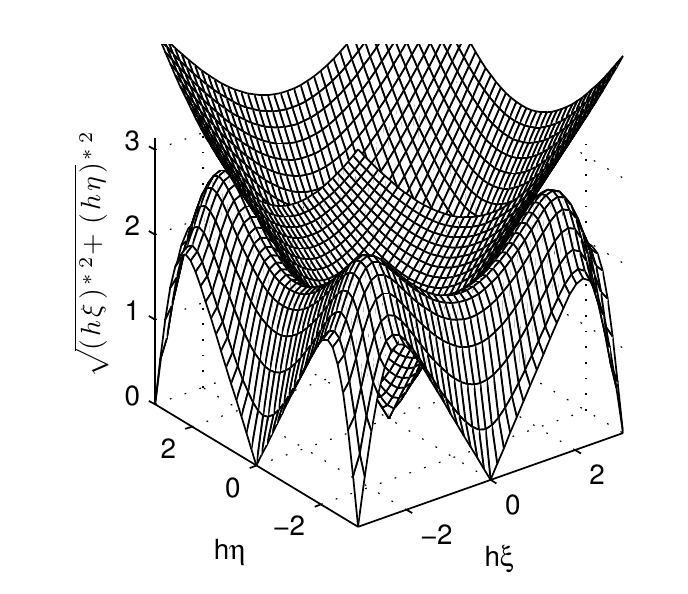}
      } \\
      a) \hspace{54mm} b) \hspace{54mm}  c)
  \end{center}
  \caption{Numerical wavenumber surfaces compared to the analytical wavenumber surface: a) second order explicit scheme (E2); b) sixth order explicit scheme (E6); c) sixth order prefactored compact scheme (Hixon). The cones represent the exact wavenumber surfaces.}
  \label{f2}
\end{figure}

In muldimensions, the numerical wavenumber and the numerical phase and group velocity are also dependent on the direction of propagation. Figure \ref{f2} shows the numerical wavenumber surface for the wave equation in two dimension, corresponding to schemes E2, E6 and Hixon as given in table \ref{t1} and equation (\ref{jjj}), respectively. The cone represents the exact wavenumber surface, obtained by revolving the straight line from figure \ref{f1} around the vertical axis. One can clearly notice the anisotropy in the numerical wavenumber surfaces associated with the finite differencing.

A simple way to reveal the numerical anisotropy is by considering the advection equation in two dimensions,

\begin{equation}\label{hh}
\partial_t u
=\textbf{c} \nabla u,
\end{equation}
with the initial condition $u(\textbf{r},0)=u_0(\textbf{r})$, where $\textbf{r}=(x,y)$ is the vector of spatial coordinates, $\textbf{c}=c(\cos\alpha \hspace{2mm} \sin\alpha)$ is the velocity vector ($c$ is a scalar and $\alpha$ the propagation direction angle), $\nabla=(\partial_x \hspace{2mm} \partial_y)^T$ and $u(\textbf{r},t)$ and $u_0(\textbf{r})$ are scalar functions. A simple semi-discretization of equation (\ref{hh}) on a square grid is obtained as

\begin{equation}\label{e2}
    d_t u=-\frac{c}{2 h} \big[\cos \alpha
(u_{{i+1},j}-u_{{i-1},j})+ \sin \alpha
(u_{i,{j+1}}-u_{i,{j-1}})\big],
\end{equation}
where $h$ is the grid step. Consider the Fourier-Laplace transform:

\begin{equation}\label{e4}
    \tilde{u}(\xi,\eta, \omega)=\frac{1}{(2\pi)^{3}}\int_{0}^{\infty} \int
\int_{-\infty}^{\infty}u(x,y,t) e^{-i (\xi x +\eta y - \omega t)} dx dy dt
\end{equation}
where $\xi = K \cos \alpha$ and $\eta = K \sin \alpha$ are the components of the
wavenumber and $\omega$ is the frequency ($K$ is the wavenumber magnitude). The application of Fourier-Laplace transform to equation
(\ref{hh}) gives the exact dispersion relation:

\begin{equation}\label{e6}
    \omega=cK(\cos^2 \alpha+\sin^2 \alpha)=cK.
\end{equation}
The exact phase velocity is given by $c_{e}=\omega/K=c$. By substituting $\omega$ in equation (\ref{e4}) with (\ref{e6}), $u(\textbf{r},t)$ is obtained as a superposition of sinusoidal solutions in the plane with constant phase lines given by $x\cos\alpha+y\sin\alpha-c_{e}t=const$. As one can notice, the exact phase velocity $c_{e}$ does not depend on the propagation direction $\alpha$, which means that the wave propagates with the same phase velocity in all directions (it is isotropic). Moreover, the exact group velocity defined as $g_{e}=\partial\omega/ \partial K=c$ is the same as the exact phase velocity because the dispersion relation is a linear function of $K$.

We now apply the same Fourier-Laplace transform to the numerical approximation  (\ref{e2}) and obtain the numerical dispersion relation in the form

\begin{equation}\label{ii}
    \omega=\frac{c}{h}\big[\cos\alpha \sin(Kh\cos\alpha)+\sin\alpha \sin(Kh\cos\alpha)\big]
\end{equation}
The numerical phase velocity will be given as

\begin{equation}\label{iii}
    c_{n}=\frac{\omega}{K}=\frac{c}{Kh}\big[\cos\alpha \sin(Kh\cos\alpha)+\sin\alpha \sin(Kh\cos\alpha)\big].
\end{equation}    
                                                                                                                                                                                                                                                                                                                                            
\begin{figure}[H]
  \begin{center}
    \mbox{
       \label{f3a}\includegraphics[width=0.4\textwidth]{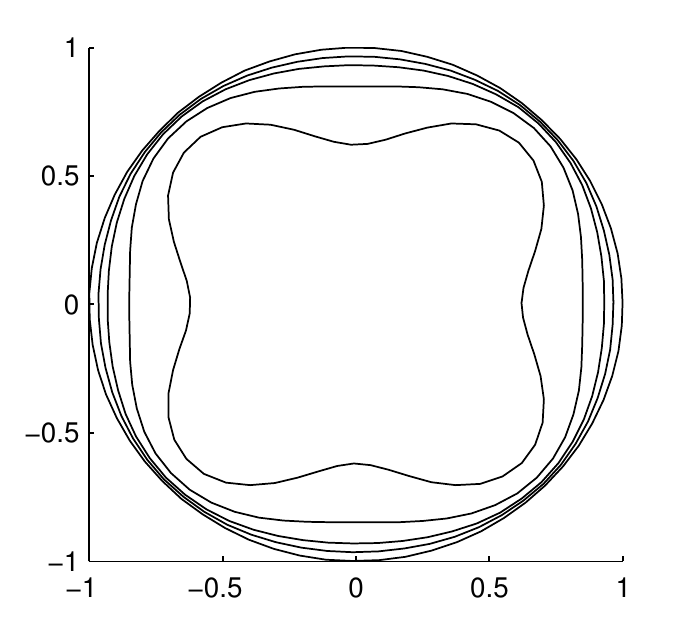}
       \label{f3b}\includegraphics[width=0.4\textwidth]{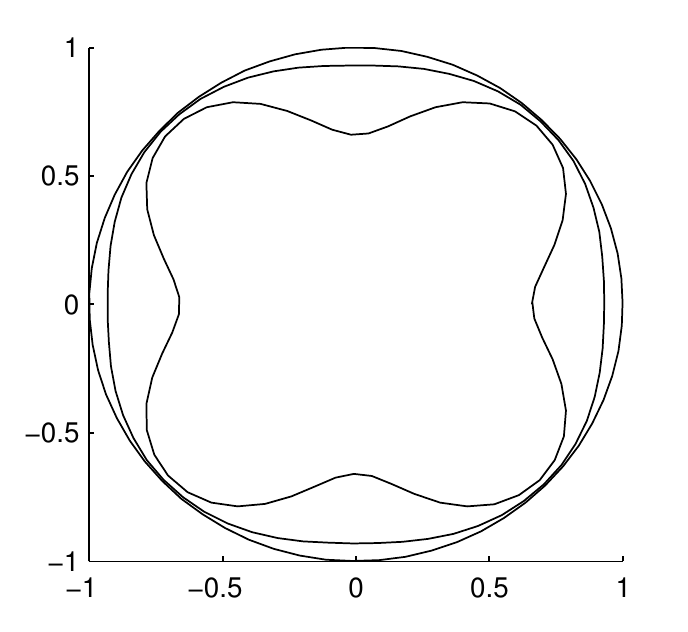}
      } \\
      a) \hspace{73mm} b)
  \end{center}
  \caption{Polar diagram of normalized phase velocities as a function of points per wavelength (PPW)
   and the direction of propagation: a) fourth-order explicit schemes (lowest number of points per wavelength is 4); b) sixth-order compact schemes (lowest number of points per wavelength is 3).}
  \label{f3}
\end{figure}

The constant phase lines are expressed by the equation $x\cos\alpha+y\sin\alpha-c_{n}t=const$ and move with the phase velocity $c_{n}$. The numerical anisotropy is revealed in equation (\ref{iii}) by the dependence of the numerical phase velocity on the propagation direction angle $\alpha$. In addition, the numerical group velocity is different from the numerical phase velocity (while previously, in the continuous case, they were the same),

\begin{equation}\label{}
    g_{n}=\partial_K \omega=c\big[\cos^2 \alpha\cos(Kh\cos\alpha)+\sin^2 \alpha\cos(Kh\sin\alpha
   )\big],
\end{equation}
which is also dependent on the propagation direction. This directional dependence of both phase and group velocities defines the numerical anisotropy. As an illustration, figure \ref{f3} shows polar diagrams for two typical schemes, fourth order explicit E4 and sixth order compact C6 schemes, revealing the numerical anisotropy (the circle of radius $1$ in figure \ref{f3} represents the exact solution). 

\section{Reduction of the Numerical Anisotropy}
\label{sec:2}

In this section, several attempts to reduce the numerical anisotropy, performed by various research groups over the years, are briefly reviewed. The optimizations of the schemes are grouped according to the mathematical model: wave equation, Helmholtz equations, advection equation, Maxwell equation, and dendritic solidification equations.

%*****************************
\subsection{Wave Equation}

Although the behavior of the numerical anisotropy was often reported in various one-dimensional optimizations of finite difference schemes, one of the first systematic attempts to specifically reduce the numerical anisotropy in finite difference schemes was introduced by Trefethen \cite{trefethen} in the framework of wave equation. To illustrate Trefethen's approach, let us consider the two dimensional wave equation in the form

\begin{equation}\label{sas}
\partial_{tt}u
= \partial_{xx}u + \partial_{yy}u,
\end{equation}
defined in $R^2 \times [0,\infty)$, with appropriate initial and boundary conditions. Using the Fourier-Laplace transform, it is ease to find the exact dispersion relation in the form $\omega^2 = \xi^2 + \eta^2$, where $\omega$ is the frequency and $(\xi, \eta)$ is the wavenumber vector. Equation (\ref{sas}) was discretized by Trefethen \cite{trefethen} on a Cartesian grid, using second order accurate schemes for both temporal and spatial derivatives as

\begin{equation}\label{}
u_{ij}^{n+1} - u_{ij}^{n} + u_{ij}^{n-1} 
= \frac{k^2}{h^2} (u_{i+1,j}^{n} + u_{i-1,j}^{n} + u_{i,j+1}^{n} + u_{i,j-1}^{n} - 4u_{i,j}^{n})
\end{equation}
which was labeled $LF^2$. Then the same scheme was used to discretize equation (\ref{sas}), except the spatial derivatives were approximated along the diagonal directions with the space step $\sqrt{2}h$; this latter discretization was termed $LF^2$. It was found that the weighted averaging $2/3 LF^2 + 1/3 LF_2$ provided a low numerical anisotropy in the order of $(\sqrt{\xi^2+\eta^2}h)^4$. Slightly the same approach was used by Vichnevetsky \cite{vich} who corrected the numerical isotropy of the wave propagation in two dimensions using either the linear advection equation or the wave equation.

In a series of papers, Sescu et al.  \cite{sescu1,sescu2,sescu3} proposed a technique to derive explicit multidimensional finite difference schemes for wave equation and Euler equations. By using the transformation matrix between two orthogonal reference frames, one aligned with the grid line and the other along the diagonal direction, the multidimensional finite difference scheme was obtained as

\begin{equation}\label{}
\left( \partial_x u \right)_{i,j}
=
\frac{1}{h(1+\beta)}
\sum_{\nu=-M}^{\nu=M} a_\nu
\left(
\textbf{E}_{x}^{\nu} +
\frac{\beta}{2}
\textbf{D}_{x}
\right)\cdot u_{i,j}
\end{equation}
where the multidimensional space shift operator $\textbf{E}_{x}^{\nu}\cdot u_{i,j} =
u_{i+\nu,j}$ (see Vichnevetsky and Bowles~\cite{vich} for one dimension) is used. The coefficients $a_n$ are those from the classical centered explicit schemes. The operator $\textbf{D}_{x}^{\nu}\cdot$  was defined as
$\textbf{D}_{x}^{\nu}\cdot =
\left(
\textbf{E}_{x}^{\nu}\textbf{E}_{y}^{\nu} +
\textbf{E}_{x}^{-\nu}\textbf{E}_{y}^{\nu}
\right)\cdot
$
The parameter $\beta$ is called isotropy corrector factor (ICF). The application of the Fourier transform to the multidimensional schemes gives the numerical wavenumber

\begin{equation}\label{}
    (\xi h)_{opt}^*=\frac{2}{(1+\beta)} \sum _{n=-N}^M a_n \Big\{e^{nI\xi h}+\frac{\beta}{2}\big[e^{nI(\xi+\eta) h}+e^{nI(\xi-\eta)
    h}\big]\Big\},
\end{equation}
Then the numerical dispersion relation corresponding to two-dimensional wave equation was considered in the form
$
    \omega^2-\big[(\xi h)_{opt}^{* \hspace{2 mm} 2} +(\eta h)_{opt}^{* \hspace{2 mm} 2}\big]=0,
$
and the ICF was determined by minimizing the integrated error between the phase or
group velocities defined along $x$ and $x=y$ directions. Two
curves in wavenumber-frequency space were considered: one was the
intersection between the numerical dispersion relation surface and
$\eta=0$ plane, and the other was the intersection between the
numerical dispersion relation surface and the $\xi=\eta$ plane. These
two curves were superposed in the $(Kh, \omega)$  plane, where $K h=\big[(\xi h)^2+(\eta h)^2\big]^{\frac{1}{2}}$.
Assuming that the equations of the two curves in $(K h,
\omega)$ plane are $\omega_1=\omega_1 (Kh,\beta)$ and $\omega_2=\omega_2 (Kh,\beta)$, the integrated error between the phase velocities was then calculated
on a specified interval as
$
    C(\beta)=\int_0^ \eta\big|c_1(Kh,\beta)-c_2(Kh,\beta)\big|^2d(Kh),
$
where $c_1(Kh,\beta)$ and $c_2(Kh,\beta)$ are the numerical phase velocities. The minimization was done by equating the first derivative of $C(\beta)$ or $G(\beta)$ with zero, which provided the value of ICF, $\beta$.

Sescu et al. \cite{sescu4,sescu5} conducted a comprehensive stability analysis of the multidimensional schemes combined with either linear-multistep or multi-stage time marching schemes, and obtained several noteworthy results. For the Leap-Frog scheme applied to the advection equations, it was shown that the stability restriction corresponding to multidimensional schemes differs from the corresponding stability restriction via conventional schemes by the factor $(2\beta+2)/(\beta+2)$, where $\beta$ is the isotropy corrector factor. The conclusion was that the stability restrictions corresponding to multidimensional schemes are more convenient compared to the conventional schemes. For an arbitrary direction of the convection velocity with $|c_x| \geq |c_y|$, the stability restriction for conventional stencils was given by $\sigma_x+\sigma_y \leq CFL$, where $\sigma_x=k|c_x|/h$ and $\sigma_y=k|c_y|/h$. For multidimensional stencils the stability restriction was given by $ (1+\beta)\sigma_x+\sigma_y \leq CFL(1+\beta)$ (where, for example, $CFL$ is $1$, $0.72874$ or $0.63052$ corresponding to E2, E4 or E6 scheme, respectively). Adams-Bashforth and Runge-Kutta time marching schemes in combination with conventional and multidimensional schemes were also analyzed, and it was found that the multidimensional schemes provide less restrictive stability limits.

%*****************************
\subsection{Helmholtz Equation}

Tam and Webb \cite{tam2} performed an anisotropy correction of the finite difference representation of the Helmholtz equation,

\begin{equation}\label{zz}
\nabla^2 p + \xi^2 p = f
\end{equation}
where $p$ is the pressure perturbation, $\nabla^2$ is the Laplacian operator, $f$ is the source distribution (e.g., a monopole), $\xi = 2\pi/\lambda$ is the wavenumber, and $\lambda$ is the acoustic wavelength. Tam and Webb \cite{tam2} showed that the finite difference discretization of the Helmholtz equation,

\begin{equation}\label{zx}
\frac{p_{i+1,j} - 2p_{i,j} + p_{i-1,j}}{h^2} + \frac{p_{i,j+1} - 2p_{i,j} + p_{i,j-1}}{h^2} 
 + \xi^2 p_{i,j} = f_{i,j}
\end{equation}
with five grid points per wavelength introduces significant numerical anisotropy (equally-spaced grid is assumed in both x- and y-directions, and the spatial step is denoted as before by $h$). They constructed an anisotropy correction factor using asymptotic solutions to the continuous equation (\ref{zz}) and its finite difference approximation (\ref{zx}) as

\begin{equation}\label{}
p_a(r,\theta)_{r_{ij}\rightarrow \infty} =
\left(
\frac{2\pi}{\xi}
\right)
\frac{\pi}{i r^{1/2}} e^{i(\xi r - \pi/4)} \bar{F}(\bar{\alpha}_s,\bar{\beta}_{+}(\bar{\alpha}_s)) + O(r^{-3/2})
\end{equation}
and

\begin{equation}\label{}
 p_n(r_{ij},\theta_{ij})_{r_{ij}\rightarrow \infty} = \frac{e^{iK_{ij}r_{ij}}}{r^{1/2}_{ij}}
\left[
G_0(\theta_{ij}+\frac{G_1(\theta_{ij}}{r_{ij}}) + O(r^{-5/2}_{ij})
\right]
\end{equation}
respectively, where $(r_{ij},\theta_{ij})$ are polar coordinates, $K_{ij} = \alpha_s(\theta_{ij})\cos{\theta_{ij}} + \beta_s(\theta_{ij})\sin{\theta_{ij}}$ (with $\alpha_s$ and $\beta_s$ being the wavenumber components from the Fourier transform), and $G_0(\theta_{ij})$ and $G_1(\theta_{ij})$ are functions depending on $\alpha_s$, $\beta_s$, $\theta$ and the Fourier transform $\bar{F}$ of the source term (for more details see equations (19) and (21) in Tam and Webb \cite{tam2}). The anisotropy corrector factor was then defined by the ratio between the absolute values of the two,

\begin{equation}\label{}
 D(\theta,\xi h) = \frac{|p_a|}{|p_n|}
\end{equation}
The correction factor is independent of the distribution of sources, meaning that it can be computed once and for all types of sources. Significant reduction of the anisotropy error was obtained.

%*****************************
\subsection{Advection Equation}

Gaitonde and Shang \cite{Gaitonde} proposed a class of high-order compact difference-based finite-volume schemes which minimized the dispersion and isotropy error functions for the range of wavenumbers of interest. The starting point was the one dimensional advection equation,

\begin{equation}\label{}
 \partial_t u + \partial_x f = 0, \hspace{4mm}
 f = cu, \hspace{4mm} c > 0
\end{equation}
which was discretized using a finite volume approach as

\begin{equation}\label{}
d_t \bar{u}_i + \bar{f}_{i+1/2} - \bar{f}_{i-1/2} = 0
 \end{equation}
 where $\bar{u}$ is the average value of $u$ inside a cell, $\bar{u} = 1/h\int_{x_{i-1/2}}^{x_{i+1/2}}udx$, and $\bar{f}$ is the flux function approximating $f$, which is dependent on the values of $\bar{u}$ from neighbor cells. The reconstruction can be done by considering a primitive function $v = \int_{0}^{x}$ which must be discretized at the cell interface. Gaitonde and Shang \cite{Gaitonde} considered a five-point compact stencil in the form
 
 \begin{equation}\label{ddd}
 \alpha v_{i-1/2} + v_{i+1/2} + \alpha v_{i+3/2}
 = b\frac{v_{i+5/2} - v_{i-3/2}}{4h} + a\frac{v_{i+3/2} - v_{i-1/2}}{2h}
 \end{equation}
 where $\alpha$, $a$, and $b$ are constants which determine the order of accuracy of the scheme. Using Taylor series expansions, they sacrificed the order of accuracy of the schemes by writing $a$ and $b$ as functions of $\alpha$,
 
 \begin{equation}\label{}
 a = \frac{2(2+\alpha)}{3}, \hspace{4mm} b = \frac{-1+4\alpha}{3}
 \end{equation}
 The spectral function associated with the scheme (\ref{ddd}) is given as
 
  \begin{equation}\label{}
 \hat{A}(w) = \frac{i\left( a\sin(w) + b \sin(2w)/2 \right)}{1 + 2\alpha \cos{w}}
 \end{equation}
 where $w = 2\pi \xi h/L$ is the scaled wave number. The dispersion error is associated with the imaginary part of the spectral function, $w_d(w) = Im(\hat{A}(w))$. A scaled isotropy wavenumber was defined as
 
\begin{equation}\label{}
w_i(w,\theta) = \cos(\theta) w_d (w \cos(\theta)) + \sin(\theta) w_d (w \sin(\theta))
\end{equation}
 where $\theta$ is the angle that the direction of propagation makes with the x-axis. An isotropy error function was defined by Gaitonde and Shang \cite{Gaitonde} in the form
 
\begin{equation}\label{}
E_i(\alpha,w_{max}) = \int_{0}^{w_{max}}\int_{0}^{\pi/2} |w_i - w| d\theta dw
\end{equation}
which was minimized to find the value of $\alpha_{opt}$ that gives the lowest numerical anisotropy. Numerical examples confirmed a considerable reduction of the isotropy error.

Sescu and Hixon \cite{sescu6,sescu7} extended the previous optimization performed in \cite{sescu2} to prefactored compact finite difference schemes \cite{hixon1,hixon2} applied to the advection equation. The prefactored compact schemes are defined on a three-point stencil and can return up to eight order of accuracy (see equations (\ref{jjj})). They can be used within a predictor-corrector type time marching scheme framework (MacCormack \cite{MacCormack}), because the numerical derivatives are determined by sweeping from one boundary to the other, in both directions. Following the same analysis as in the case of explicit schemes, the multidimensional prefactored compact schemes were obtained as

\begin{eqnarray}\label{e1}
u_{i,j}^{F'} &=& \frac{\alpha}{1+\beta} \left[ u_{i+1,j}^{F'} + \frac{\beta}{2} \left( u_{i+1,j-1}^{F'} + u_{i+1,j+1}^{F'} \right) \right]  \\
&+& \frac{1}{h(1+\beta)}
\left[
b u_{i+1,j} - e u_{i,j} 
+ \frac{\beta}{2}
\left(
b u_{i+1,j+1} + b u_{i+1,j-1} - 2e u_{i,j} 
\right)
\right] \nonumber \\
u_{i,j}^{B'} &=& \frac{\alpha}{1+\beta} \left[ u_{i-1,j}^{B'} + \frac{\beta}{2} \left( u_{i-1,j-1}^{B'} + u_{i-1,j+1}^{B'} \right) \right]  \\
&+& \frac{1}{h(1+\beta)}
\left[
e u_{i,j} - b u_{i-1,j}
+ \frac{\beta}{2}
\left(
2e u_{i,j} - b u_{i-1,j+1} - b u_{i-1,j-1}
\right)
\right] \nonumber
\end{eqnarray}
for fourth order of accuracy, and

\begin{eqnarray}\label{e1}
u_{i,j}^{F'} &=& \frac{\alpha}{1+\beta} \left[ u_{i+1,j}^{F'} + \frac{\beta}{2} \left( u_{i+1,j-1}^{F'} + u_{i+1,j+1}^{F'} \right) \right]  \\
&+& \frac{1}{h(1+\beta)}
\left[
b u_{i+1,j} - e u_{i,j} - f u_{i-1,j} 
+ \frac{\beta}{2}
\left(
b u_{i+1,j+1} - f u_{i-1,j-1} +
b u_{i+1,j-1} - f u_{i-1,j+1} - 2e u_{i,j}
\right)
\right] \nonumber \\
u_{i,j}^{B'} &=& \frac{\alpha}{1+\beta} \left[ u_{i-1,j}^{B'} + \frac{\beta}{2} \left( u_{i-1,j-1}^{B'} + u_{i-1,j+1}^{B'} \right) \right]  \\
&+& \frac{1}{h(1+\beta)}
\left[
b u_{i+1,j} - e u_{i,j} - b u_{i-1,j}
+ \frac{\beta}{2}
\left(
f u_{i+1,j+1} - b u_{i-1,j-1} +
f u_{i+1,j-1} - b u_{i-1,j+1} - 2e u_{i,j}
\right)
\right] \nonumber
\end{eqnarray}
for sixth order of accuracy. $\beta$ is the isotropy corrector factor (ICF) and its magnitude can be determined by minimizing the dispersion error corresponding to the wave-front propagating along a grid line and the wave-front propagating along a diagonal direction.

Using Fourier analysis, the numerical wavenumbers and the numerical dispersion relation corresponding to the two dimensional wave equation were found. The individual (forward or backward) numerical wavenumber has both real and imaginary parts: the real part of the forward operator is equal to the real part of the backward operator, and the imaginary parts are opposite. As a result, in a MacCormack predictor-corrector scheme the overall imaginary part is zero. The real parts of the numerical wavenumbers corresponding to multidimensional schemes, for derivatives along $x$-direction, were given by:

\begin{eqnarray}\label{e1}
Re[(kh)^*_{m}] =
\frac{1}{1+\beta}
\left\{
f_m(\eta_x) +
\frac{\beta}{2}
\left[
f_m(\eta_x + \eta_y) +
f_m(\eta_x - \eta_y)
\right]
\right\},
\end{eqnarray}
where $m=4$ for fourth and $m=6$ for sixth order of accuracy, $f_4(\eta_x) = 3 \sin{\eta_x}/(2+\cos{\eta_x})$, $f_6(\eta_x) =
(28 \sin{\eta_x} + \sin{2\eta_x})/(18+12\cos{\eta_x})$, $\eta_x = \xi h$, $\eta_y = \eta h$ and $\xi$ and $\eta$ are the components of the wavenumber.

In terms of numerical stability, more efficient stability restrictions were obtained as in the case of multidimensional explicit schemes. For example, multidimensional MacCormack schemes were found to provide a stability restriction in the form

\begin{eqnarray}\label{r1}
[\sigma_x(1+\beta)]^{2/3} + \sigma_y^{2/3} \leq \frac{(1+\beta)^{2/3}}{\xi_{max}},
\end{eqnarray}
if $|c_x| \geq |c_y|$, and 

\begin{equation}\label{r2}
\sigma_x^{2/3} + [\sigma_y(1+\beta)]^{2/3} \leq \frac{(1+\beta)^{2/3}}{\xi_{max}},
\end{equation}
if $|c_y| \geq |c_x|$. For diagonal directions, with respect to the grid, ($|c_x|=|c_y|=|c|$) the stability restriction becomes

\begin{equation}\label{s2}
\sigma \leq 
\frac{(1+\beta)}{\xi_{max}^{3/2} \left[ 1+(1+\beta)^{2/3} \right]^{3/2}}.
\end{equation}
It is obvious that the right hand side of equation (\ref{s2}) is greater than $1/(2\xi_{max})^{3/2}$ when $\beta > 0$, and goes to $1/(\xi_{max})^{3/2}$ when $\beta \rightarrow \infty$. This generated more efficient stability restrictions by using multidimensional compact schemes. Test cases showed that the multidimensional compact schemes were more efficient for both fourth and sixth order accurate schemes.

%*****************************
\subsection{Maxwell Equations}

Sun and Trueman \cite{Sun1} performed an optimization of finite difference schemes applied to Maxwell equations, in terms of reducing the dispersion and isotropy errors. For brevity, we show here the numerical dispersion relations (for finite differencing representations of the Maxwell equations, see equations (1), (2) and (4) in Sun and Trueman \cite{Sun1}):

\begin{equation}\label{q1}
\left( \frac{\sin(\omega k/2)}{ck} \right)^2
= \left( w\frac{\sin(\beta_a k/2)}{h} + (1-w)\frac{\sin(3\beta_a k/2)}{3h} \right)^2
\end{equation}
corresponding to a grid line, and

\begin{equation}\label{q2}
\left( \frac{\sin(\omega k/2)}{ck} \right)^2
= 2\left( w\frac{\sin(\beta_d k/2)}{h} + (1-w)\frac{\sin(3\beta_d k/2)}{3h} \right)^2
\end{equation}
corresponding to the diagonal direction, where $w$ is a weighting factor, $\beta_a$ is the numerical phase constant along the grid line, $\beta_d$ is the numerical phase constant along the diagonal direction, $\omega$ is the frequency, and $k$ is the time step (an equally-spaced grid is considered again). The optimization in terms of reducing the numerical anisotropy was done by eliminating the time step terms in equations (\ref{q1}) and (\ref{q2}) to obtain

\begin{equation}\label{}
w_i = \frac
{\sqrt{2}  \sin(3\beta_d k/2)/(3h) - \sin(3\beta_a k/2)/(3h)}
{\left[ \sin(\beta_a k/2)/h - \sin(3\beta_a k/2)/(3h) \right]  
-  \sqrt{2}  \left[ \sin(\beta_d k/2)/h - \sin(3\beta_d k/2)/(3h) \right]}
\end{equation}
This optimal weight $w_i$ is a function of mesh density only, and is not dependent of the time step size or the frequency of the signal. This method theoretically provides a uniform phase velocity in all directions. Further optimizations of this scheme were performed in another paper of Sun and Trueman \cite{Sun2}.

Koh et al. \cite{Koh} derived a two dimensional finite-difference time-domain method, discretizing the Maxwell equations, to eliminate the numerical dispersion and anisotropy. The proposed scheme is given as

\begin{eqnarray}\label{ff}
d_t^2 H_{x,i,j+1/2}^n = -\frac{k}{\mu h} d_y E_{x,i,j+1/2}^n  \nonumber  \\
d_t^2 H_{y,i+1/2,j}^n = -\frac{k}{\mu h} d_x E_{y,i+1/2,j}^n   \\
d_t^2 E_{z,i,j}^{n+1/2} + \frac{\sigma k}{2\epsilon} [E_{z,i,j}^{n+1} + E_{z,i,j}^{n}]
 = \frac{k}{\epsilon h} d_x H_{y,i,j}^{n+1/2} -  \frac{k}{\epsilon h} d_y H_{x,i,j}^{n+1/2} \nonumber
\end{eqnarray}
where $d_t^2$ is the central difference operator with respect to time,

\begin{eqnarray}\label{}
d_p f_q = \left( 1-\frac{\alpha}{2} \right) d_p f_q + \frac{\alpha}{4} \left( d^2_p f_{q+1} + d^2_p f_{q-1} \right)
\end{eqnarray}
with $p$ or $q$ being either $x$ or $y$, and

\begin{eqnarray}\label{}
d^2_x f_{i,j} = f_{i+1/2,j} - f_{i-1/2,j}, \hspace{4mm} d^2_y f_{i,j} = f_{i,j+1/2} - f_{i,j-1/2}
\end{eqnarray}
where $f$ is a generic function. In equation (\ref{ff}), $E$ is the electric field, $H$ is the magnetic field strength, $\sigma$, $\mu$ and $\epsilon$ are the conductivity, permeability and the permittivity, respectively, of the domain, $k$ is the time step, and $h$ is the spatial step in all directions. For a nonconductive media $\sigma = 0$, the numerical dispersion relation of can be obtained as

\begin{eqnarray}\label{mm}
\frac{1}{h^2} C_{+} C_{\times}  \left( \alpha - \frac{2}{C_{+}} \right)^2
- \frac{1}{h^2} \left( \frac{4C_{\times}}{C_{+}} - C_{+} \right)
- \frac{1}{(ck)^2} \sin^2 \left( \frac{\omega k}{2} \right)
\end{eqnarray}
where $C_{+} = \sin^2 (\xi h/2) + \sin^2 (\eta h/2)$, $C_{\times} = \sin^2 (\xi h/2)  \sin^2 (\eta h/2)$, and $\xi$ and $\eta$ are the components of the wavenumber. Equation (\ref{mm}) is a quadratic equation in $\alpha$, and the solution is given as

\begin{eqnarray}\label{mm}
\alpha = \frac{2}{C_{+}}
\left[
1 - \sqrt{1 - \frac{h^2 C_{+}}{4C_{\times}} \left( \frac{1}{h^2} C_{+} - \frac{1}{(ck)^2} \sin^2 \left( \frac{\omega k}{2} \right) \right)}
\right]
\end{eqnarray}
An optimal value for $\alpha$, achieving an isotropic numerical phase velocity, can be simply estimated as the mean value of $\alpha$ over the azimuthal angles, and it was found that it remains constant (approximately, $0.167$) for a wide range of grid sizes, and it is insensitive to the value of the Courant number.

Kim et al. \cite{Kim} derived new three-dimensional isotropic dispersion-finite-difference time-domain schemes (ID-FDTD) based on a a linear combination of the traditional central difference equation and a new difference equation based on the extra sampling points. They used the same scaling factors as for the two-dimensional case to attain isotropic dispersion and exact phase velocity. Based on the weighting factors, seven different FDTD schemes were formulated, including the Yee scheme \cite{Yee}. Among the seven proposed FDTD schemes, three showed improved isotropy of the dispersion compared to the dispersion of the Yee scheme. For the sake of brevity, the complete expressions of the schemes are not included here (see Kim et al. \cite{Kim} for more details), and only the numerical dispersion relation is briefly presented. Plane wave solutions were introduced in discretized forms as

\begin{eqnarray}\label{ii1}
\textbf{E}^n_{i,j} = \textbf{E}_0 e^{I(n\omega k - \xi i h - \eta j h - \zeta k h)} 
\end{eqnarray}

\begin{eqnarray}\label{ii2}
\textbf{H}^n_{i,j} = \textbf{H}_0 e^{I(n\omega k - \xi i h - \eta j h - \zeta k h)}
\end{eqnarray}
where $I = \sqrt{-1}$, $\omega$ is the frequency, $(\xi,\eta,\zeta)$ is the numerical wavenumber vector, and $\textbf{E}_0$ and $\textbf{H}_0$ are constant vectors. After inserting (\ref{ii1}) and (\ref{ii2}) into the discretized form of the Maxwell equations (see equation (10) in Kim et al. \cite{Kim}), matrix equations are obtained as
$
C \textbf{H}_0 = S_t \epsilon_0 \textbf{E}_0, \hspace{1mm}
C \textbf{E}_0 = S_i \mu_0 \textbf{H}_0
$
where

\begin{eqnarray}
C = \left[
\begin{array}{ccccc}
 0      &   -K_z   &   K_y \\
 K_z  &    0       &  -K_x \\
-K_y  &    K_x   &   0
\end{array}
\right]
\end{eqnarray}
and $K_p = S_p/h [\alpha (P_p - Q_p) - \beta Q_p/2 + 1]$ ($p$ being either $x$, $y$ or $z$), $S_x = \sin(\xi h/2)$, $S_y = \sin(\eta h/2)$, $S_z = \sin(\zeta h/2)$, $P_x = Sy Sz$, $P_y = Sx Sz$, $P_z = Sx Sy$, $Q_x = S_y^2 + S_z^2$, $Q_y = S_x^2 + S_z^2$, $Q_z = S_x^2 + S_y^2$, and $S_t = \sin{\omega k/2}/k$. An eigenvalue equation was obtain as

\begin{eqnarray}\label{}
(C^2 + S_t^2 \mu_0 \epsilon_0 I) = 0,
\end{eqnarray}
and the numerical dispersion relation was obtained by vanishing the associated determinant,

\begin{eqnarray}\label{}
\frac{S_t^2}{c_0^2} = K_x^2 + K_y^2+ K_z^2
\end{eqnarray}
where $c_0 = 1/\sqrt{\epsilon_0 \mu_0}$. The isotropy correction was performed by defining the values of the weighting factors $\alpha$ and $\beta$, which unlike the two-dimensional case are not unique. Kim et al. \cite{Kim} used the scaling factor from the two-dimensional case, and modified the numerical dispersion relation to estimate the weighting factors.

%*****************************
\subsection{Dendritic Solidification}

Kumar \cite{Kumar} derived isotropic finite difference schemes for the first and second derivatives in the context of symmetric dendritic solidification. The first derivative was discretized as

\begin{eqnarray}\label{vv}
(\partial_x u)_{I,i,j} = \frac{1}{2h}
\left[
\frac{1}{6}(u_{i+1,j+1} - u_{i-1,j+1})
\right.   \nonumber \\
\left.
 + \frac{4}{6}(u_{i+1,j} - u_{i-1,j})
\right. \\
\left.
 + \frac{1}{6}(u_{i+1,j-1} - u_{i-1,j-1})
\right]   \nonumber
\end{eqnarray}
which involves grid points not only along $x$-direction, but also along $y$-direction. The Taylor expansion of the scheme (\ref{vv}) can be written as $(\partial_x u)_{I,i,j} = (1+ h^2/6 \nabla^2)(\partial_x u)_{i,j}$, where the leading order term involves the Laplacian only, implying no directional dependence. The second derivative was discretized as

\begin{eqnarray}\label{}
(\partial_{xx} u)_{I,i,j} = \frac{1}{h^2}
\left[
\frac{1}{12}(u_{i+1,j+1} - 2u_{i,j+1} u_{i-1,j+1})
\right.   \nonumber \\
\left.
+ \frac{10}{12}(u_{i+1,j} - 2u_{i,j} + u_{i-1,j})
\right.  \\
\left.
+ \frac{1}{12}(u_{i+1,j-1} - 2u_{i,j-1} + u_{i-1,j-1})
\right]   \nonumber
\end{eqnarray}
where the Taylor expansion is given by $(\partial_{xx} u)_{I,i,j} = (1+ h^2/12 \nabla^2)(\partial_{xx} u)_{i,j}$, being again a function of the Laplacian only. The conventional cross derivative $(\partial_{xy} u)_{I,i,j}$ was found to be intrinsically isotropic according to the criterion developed by Kumar \cite{Kumar}. The Laplacian can be obtained by combining the isotropic derivatives along x- and y-directions, $(\nabla^2 u)_{i,j} = (\partial_{xx} u)_{I,i,j} + (\partial_{yy} u)_{I,i,j}$. Significant reduction of the numerical anisotropy was obtained by using these schemes. Shen and Cangellaris \cite{Shen} exploited further this approach to develop new isotropic finite-difference time-domain schemes modeling electromagnetic wave propagation.

%\section{Future Directions}
%\label{sec:4}

\section{Concluding Remarks}

Numerical anisotropy in finite difference discretizations of partial differential equations was discussed and reviewed. In some instances, the numerical anisotropy can be neglected, and the focus is directed toward other types of one-dimensional errors, such as numerical dispersion, dissipation or aliasing. These errors can be analyzed in the context of one dimensional differencec equations, while the extension to multidimensional discretizations is straightforward. By increasing the accuracy of one dimensional schemes or by increasing the number of grid points in the grid, the isotropic characteristics of the waves in multi-dimensions can be improved. These two practices, however, are not always effective since an increase in accuracy may require larger stencils which may introduce spurious waves at the boundaries of the domain, while by increasing of the resolution of the grid may increase the computational time. It is necessary then to analyzed the schemes in multi-dimensions and design specific optimizations with the specific objective of reducing the numerical anisotropy, and at the same time of conserving the dispersion characteristics of the corresponding one dimensional schemes. Various attempts to reduce the numerical anisotropy in finite differencing applied to various model equations were presented and discussed.

Future directions should focus on optimizations of existing compact finite difference schemes in terms of reducing the numerical anisotropy, or derivations of novel compact schemes with low numerical anisotropy. Optimizations and derivations of finite volume schemes (in terms of reducing the numerical anisotropy) applied to either structured or unstructured grids should be also taken into account, especially in the framework of wave propagation problems. Filtering schemes, as applied, for example, in large eddy simulations to separate the small scales from the large scales, may experience numerical anisotropy since they are effective at high wavenumber ranges. Optimizations of such filters in terms of reducing the numerical anisotropy is also another future area of research.

%\section{Acknowledgments}

%\begin{acknowledgements}
%\end{acknowledgements}



\begin{table}[htpb]
 \begin{center}
  \caption{Weights of the selected spatial finite difference stencils}
  \label{t1}
  \begin{tabular}{rrrrrrrrrr} \hline
       Stencil & $\alpha_1$ & $\alpha_2$ & $a_1$ & $a_2$ & $a_3$   \\\hline
       $E2$ &  0 & 0 &  1/2 & 0 & 0 \\
       $E4$ &  0 &0 &  2/3 & -1/12 & 0 \\
       $E6$ &  0& 0 &  3/4 & -3/20 & 1/60 \\
       $DRP$ &  0 &  0 &  0.770882380 &  -0.166705904 & 0.020843142 \\
       $C4$ &  1/4& 0 &  3/4 & 0 & 0 \\
       $Haras$ &  0.3534620& 0 &  1.5669657/2 & 0.13995831/4 & 0 \\
       $Lui$ & 0.5381301& 0.0666331 &  1.36757772/2 & 0.823428170/4 & 0.0185207834/6\\
       $Lele$ &  0.5771439& 0.0896406 &  1.3025166/2 & 0.99355/4 & 0.03750245/6 \\
\hline
  \end{tabular}
 \end{center}
\end{table}

\end{document}